\documentstyle[leqno]{article}

\input{amssym}
\input{arrow}
\def\di{\displaystyle }
\def\be{\begin{eqnarray}}
\def\ee{\end{eqnarray}}
\begin{document}
\title{Self-equivalence 3rd order ODEs by time-fixed transformations}
\author{ Mehdi Nadjafikhah \and Ahmad Reza Forough }
\date{}
\maketitle
\begin{abstract}
Let $y'''=f(x,y,y',y'')$ be a 3rd order ODE. By Cartan equivalence
method, we will study the local equivalence problem under the
transformations group of time-fixed coordinates.

\medskip \noindent {\it A.M.S. 2000 Classification Number:} 58A15
\end{abstract}
\section{Introduction}
Cartan's method of equivalence (see \cite{1}, \cite{2} and
\cite{3}) is acknowledged to be a powerful tool for studying
differential invariants. The main goal of this method is to find
necessary and sufficient conditions in order that two geometric
structures be equivalent, by a class of given diffeomorphisms. By
introducing the invariance of the differential equation under a
continuous group of symmetries, {\it Sophus Lie} rose to the
challenge of finding a general method to uncover such invariants,
but his approach had some serious defects. Roughly speaking a
symmetry group of a system of differential equations is a group
which transforms solutions of the system to another solutions. In
the classical framework of Lie, these groups consist of
geometrical transformations on the space of independent and
dependent variables for the system, and act on the solutions by
transforming their graphs. Constructing the compatible coframes
were the main part of this method, and was done by {\it \'Elie
Cartan}.

In the first step of his attempt, E. Cartan introduced the
structure equations, which leads him to the  differential
invariants. In this paper, we study Cartan's equivalence problem
$y'''=f(x,y,y',y'')$ under the transformations group \be
X=x,\;\;\;Y=\varphi(y).\label{eq:trans} \ee This is called {\it
time-fixed geometry of} the 3rd order ODEs.
\section{Cartan's equivalence problem}
Let $\omega_U=(\omega^i_U)$ and $\Omega_V=(\Omega^i_V)$ be two
coframes on an open sets $V\subset{\Bbb R}^n$ and $U\subset{\Bbb
R}^n$ respectively, and $G\subseteq {\rm GL}(n;{\Bbb R})$ be a
prescribed linear group, then find necessary and sufficient
conditions that there exist a diffeomorphism $\Phi:U\rightarrow V$
such that for each $u\in U$ \be \Phi^*\Omega_V\Big|_{\Phi(u)} =
\gamma_{VU}(u).\omega_U\Big|_u, \ee where $\gamma:U\rightarrow G$.
(In the future we will always omit the base point notation and
write the last relation as $\Phi ^* \Omega_V =
\gamma_{VU}.\omega_U$.)
\section{Time-fixed problem}
Let $(U,x,y,y',y'')$ and $(V,X,Y,Y',Y'')$ be open sets with
standard coordinates on the $2-$jet bundle $J^2({\Bbb R};{\Bbb
R})$ of mappings ${\Bbb R}\to {\Bbb R}$, and let there be given
3rd order ODEs
\begin{eqnarray}
y'''=f(x,y,y',y''),\hspace{1cm} \mbox{and} \hspace{1cm}
Y'''=F(X,Y,Y',Y''). \end{eqnarray}

The usual symmetries of these equations are the diffeomorphisms
$\Phi(x,y,y',y'')$ which map the integral curves into integral
curves, that is
\begin{eqnarray}
\Phi^*\left( \begin{array}{c} dY-Y'\,dX \\
dY'-Y''\,dX \\ dY''-F\,dX \end{array} \right) = \left(
\begin{array}{ccc}  m & 0 & 0
\\ n & p & 0 \\ q & r & s \end{array} \right) \left( \begin{array}{c} dy-y'\,dx \\ dy'-y''\,dx
\\dy''-f\,dx\end{array} \right),
\end{eqnarray}
where $mps\neq0$. By the transformation (\ref{eq:trans}), we also
have the following Jacobian condition on the diffeomorphisms:
\begin{eqnarray}
\Phi^*\left( \begin{array}{c} dX \\ dY \end{array} \right) =
\left( \begin{array}{cc} 1 & 0 \\ 0 & u \end{array} \right) \left(
\begin{array}{c} dx \\ dy \end{array} \right),
\end{eqnarray}
where $u=\varphi'(y)$. Since $\dim J^2({\Bbb R};{\Bbb R})=4$ and
there exist five relations, this is an over-determined problem on
the generators: \be dX,\;\; dY ,\;\; dY-Y'\,dX ,\;\;
dY'-Y''\,dX,\;\;dY''-F\,dX  \ee
\section{Solving the problem}
The relation \be (dY-Y'\,dX)-dY+Y'\,dX=0,\ee would seem to suggest
modifying the forms, hence define: \begin{eqnarray}
&& \Omega^1_V:=dX,\;\; \Omega^2_V:=\frac{dY}{Y'},\;\;\Omega^3_V:=\frac{dY'-Y''\,dX}{Y'}, \\
&& \Omega^4_V:=dY''-F\,dX,\;\;
\Omega^5_V:=\frac{dY-Y'\,dX}{Y'}\nonumber
\end{eqnarray}

Now we have the following relation between forms: \be
\Omega^4_V-\Omega^2_V+\Omega^1_V=0.\ee By the Jacobian conditions,
we can find the diffeomorphisms which satisfy in the above ones;
now by (\ref{eq:trans}), we have
\begin{eqnarray}
\Phi^*\left( \begin{array}{c} \Omega^1_V \\[1mm] \Omega^2_V\\[1mm]
\Omega^3_V\\[1mm]\Omega^4_V \end{array} \right) = \left( \begin{array}{cccc}
a_{11} & 0 & 0 & 0 \\ a_{21} & a_{22} & 0 & 0 \\ a_{31} & a_{32} &
a_{33} & 0 \\ a_{41} & a_{42} & a_{43} & a_{44}
\end{array} \right) \left( \begin{array}{c} \omega^1_U \\[1mm] \omega^2_U \\[1mm]
\omega^3_U \\[1mm] \omega^4_U \end{array} \right).
\end{eqnarray}
Then, $a_{11}=1$ since
\begin{eqnarray}
\Phi^*\Omega^1_V = \Phi^*dX= dx = \omega^1_U,
\end{eqnarray}
in the same manner, $a_{21}=0$ and $a_{22}=1$, since
\begin{eqnarray}
\Phi^*\Omega^2_V &=& \Phi^*\left(\frac{dY}{Y'}\right) = \frac{d(\varphi(y))}{((\varphi(y))_x} \nonumber \\
&=& \frac{\varphi'(y)\,dy}{\varphi'(y)\,y'} =\omega^2_U
\end{eqnarray}
Moreover, $a_{32}=-a_{31}=\varphi'(y)=u$, because
\begin{eqnarray}
\Phi^*(\Omega^3_V) &=&\Phi^*\left(\frac{dY'-Y''\,dX}{Y'}\right)=\frac{\varphi'(y)\,dy'+y'\,\varphi''(y)\,dy}{\varphi'(y).y'}\nonumber\\
&=&\frac{dy'-y''\,dx}{y'}+\frac{y'\,\varphi''(y)}{\varphi'(y)}.\frac{dy}{y'}-\frac{y'\,\varphi''(y)}{\varphi'(y)}.dx\\
&=& \omega^3_U+u.\omega^2_U-u.\omega^1_U. \nonumber
\end{eqnarray}
Now by assuming $v=\ell y'$, we have
\begin{eqnarray}
\Phi^*\Omega^4_V &=&\Phi^*\,(dY''-F\,dX) \nonumber \\
&=& \ell.(dy-y'\,dx)+a.(dy'-y''\,dx)+b.(dy''-f\,dx) \\
&=& -v.\omega^1_U+v.\omega^2_U+a.\omega^3_U+b.\omega^4_U.
\nonumber
\end{eqnarray}
thus, $a_{42}=-a_{41}=v$, $a_{43}=a$, $a_{44}=b$; so, the group
structure $G\subset{\rm GL}(4,{\Bbb R})$ is the set of elements in
the following form:
\begin{eqnarray}
g(a,b,u,v):= \left( \begin{array}{cccc} 1 & 0 & 0 & 0 \\ 0 & 1 & 0 & 0 \\
-u & u & 1 & 0 \\ -v & v & a & b\end{array} \right),\;\;\;
(a,b,u,v\in{\Bbb R}\;,\;b\neq0)
\end{eqnarray}
\paragraph{Theorem 1.} {\it $G$ is a  $4-$dimensional Lie subgroup
of ${\rm GL}(4,{\Bbb R})$ with multiplication: \be
g(a,b,u,v).g(a',b',u',v')=g(a+ba',bb',u+u',v+au'+bv') \ee and
inversion, \be
g(a,b,u,v)^{-1}=g\big(-\frac{a}{b},\frac{1}{b},-u,\frac{a\,u-v}{b}\big)\ee
and its Lie algebra is the set of all matrices in the form:
\begin{eqnarray}
\left(
\begin{array}{cccc}
0 & 0 & 0 & 0 \\
0 & 0 & 0 & 0 \\
-u & u & 0 & 0 \\
-v & v & a & b
\end{array}
\right)\in \;{\rm Mat}(4\times 4;{\Bbb R})
\end{eqnarray}}

\medskip \noindent {\it Proof:} The first two statements are
trivial. The last part is due to defining relations on the
Maurer-Cartan matrix form and in fact the defining relations on
the Lie algebra of $G$, so it is necessary to compute $dg.g^{-1}$.
 \hfill\ $\Box$
\section{Prolongation}
Now we lift the problem to the associated spaces $U\times G$ and
$V\times G$ with the natural left action, that is \be g\cdot
(p,h)=(p,gh), \hspace{1cm} g,h\in G,\;\;p\in U\;\;
\mbox{or}\;\;p\in V.\ee Given $\Omega_V=(\Omega^i_V)$ and
$\omega_U=(\omega^i_U)$ are adapted coframes on open sets
$V,U\subseteq \Bbb R^n$ respectively, and diffeomorphism
$\Phi:U\rightarrow V$ satisfying
\begin{eqnarray}
\Phi^*\Omega_V=\gamma_{VU}.\omega_U, \hspace{1cm}
\gamma_{VU}:U\rightarrow G. \label{eq:3}
\end{eqnarray}
We define new column vectors of $1-$forms on $V\times G$ and
$U\times G$ by \be
\Omega|_{(V,g)}=g\,\Pi^*_V\,\Omega_V,\hspace{1cm}\omega|_{(U,h)}=h\,\Pi^*_U\,\omega_U
\ee respectively, where  $\Pi_V:V\times G\rightarrow V$ and
$\Pi_U:U\times G\rightarrow U$ are natural projections.
\paragraph{Theorem 2.} {\it There exists a diffeomorphism $\Phi:U\rightarrow
V$ satisfying (\ref{eq:3}) if and only if there exists a
diffeomorphism $\Phi^1:U\times G\rightarrow V\times G$ such that
$\Phi^{1*}\,\Omega=\omega$.} (See \cite{3}.)

\medskip The above theorem is the key to the usefulness of the lifting
procedure. Moreover this diffeomorphism $\Phi^1$ covers $\Phi$,
i.e. the diagram with the natural projections \be \commdiag{
U\times G            & \mapright^{\Phi^1} & V\times G \cr
\mapdown\lft{\Pi_U}  &                    & \mapdown\rt{\Pi_V}
\cr U                    & \mapright^{\Phi} & V } \ee commutes.
Further, $\Phi^1$ is uniquely determined and automatically
satisfies \be\Phi^1(u,gh)=g.\Phi^1(u,h),\hspace{1cm}g,h\in
G,\;\;u\in U.\ee
\paragraph{Definition.} $\!\!\!$({\it Right invariant Maurer-Cartan
$1-$forms.})
Assume $G$ be a Lie group and let $R_c$ denote right
multiplication by $c\in G$, if we choose a basis
$\{\omega^i\big|_e\}$ of $T^*_eG$, the cotangent space of $G$  at
the identity point $e\in G$, then we may define global
differential forms by \be
\omega^i|_A\;=\;R^*_{A^{-1}}(w^i|_e),\hspace{1cm}\forall A\in G.
\ee Since
\begin{eqnarray}
R^*_C\,(\omega^i|_{AC}) &=& R^*_C\circ R^*_{(AC)^-1}\,(\omega^i|_e) \nonumber \\
&=& R^*_C\circ R^*_{{C}^-1}\circ R^*_{{A}^-1}\,(\omega^i|_e) \\
&=& \omega^i|_A \nonumber
\end{eqnarray}
there are a basis for the right invariant Maurer-Cartan $1-$forms
. \medskip Matters being so, a set of right invariant
Maurer-Cartan $1-$forms $\{\omega^i|_e\}$ defines functions
$C^i_{jk}$ via the equations \be dw^i=\frac{1}{2}\sum_{j<k}
C^i_{jk}.\omega^j\wedge\omega^k \ee where $C^i_{jk}=-C^i_{kj}$.

The right translational invariance immediately implies that the
functions $C^i_{jk}$ are in fact constants. These constants are
called the {\it structure constants of} $G$, relative to the
choice of Maurer-Cartan $1-$forms.
\section{Absorption first step }
Define
\begin{eqnarray}
\omega^1_U &=& dx  \nonumber \\
\omega^2_U &=& \frac{dy}{y'}\\
\omega^3_U &=&
-u\,dx+u\,\frac{dy}{y'}+\big(\frac{dy'-y''\,dx}{y'}\big)\nonumber\\
\omega^4_U
&=&-v\,dx+v\frac{dy}{y'}+a\,\big(\frac{dy'-y''\,dx}{y'}\big)+b\,(dy''-f\,dx)\nonumber.
\end{eqnarray}
we drop the index $U$ and differentiate the $(\omega^i)$'s, giving
\begin{eqnarray}
d\omega^i=\sum A^i_{jk}.\Pi^k \wedge \omega^j+\sum
T^i_{jk}(u,g).\omega^j\wedge \omega^k \label{eq:4}
\end{eqnarray}
and we called them {\it structure equations}. The matrix
$A^i_{jk}\Pi^k$ are now Lie algebra valued differential form. The
terms involving the coefficients $T^i_{jk}$ are called {\it
torsion terms}, and the coefficients themselves are called the
{\it torsion coefficients}.

Equations (\ref{eq:4}) do not define the torsion coefficients nor
$1-$forms $\Pi^k$ uniquely, so it is necessary to simplify, even
eliminate if possible, this process is called {\it Lie algebra
valued compatible absorption}. So we have \be d\omega^1 = d\,(dx)
= 0,\;\;\;\;\; d\omega^2 = \sum_{i<j}
T^2_{ij}\;\omega^i\wedge\omega^j \ee where \be
T^2_{12}=-u-\frac{y''}{y'},\;\;\;\;T^2_{23}=1,\ee and  the rest
are zero. With respect to this reality that,the elements of the
group, which are only in the last two rows, are essential torsion
coefficients and thus they are used for reducing the parameters.
Now, by absorption of u, we will have $u=-y''/y'$, and the
structure group, reduced to the following subgroup,
\begin{eqnarray}
g:=\left(\begin{array}{cccc}
1 & 0 & 0 & 0 \\
0 & 1 & 0 & 0 \\
0 & 0 & 1 & 0 \\
-v & v & a & b
\end{array}\right)
\end{eqnarray}

We substitute the acquired value for group parameter $u$, and
repeat the procedure in the sequel, so we have
\begin{eqnarray}
\omega^1 &=& dx,\;\;\;\;\;\;\;\;\;\;\;\;\; \omega^2
 =\frac{dy}{y'},\;\;\;\; \nonumber\\ \omega^3 &=&
-\frac{y''\,dy}{y'^2}+\frac{dy'}{y'},\;\;\;\; \\
\omega^4 &=& \big(\frac{a\,y''}{y'}-v-b\,f \big)dx
+\frac{v\,dy}{y'}+\frac{a\,dy'}{y'}+b\,dy''\nonumber.
\end{eqnarray}
\section{Absorption second step}
By the last step, and using the group parameter $u$ ,we can
compute the torsion coefficients, in this step some parameters
eliminate,
\begin{eqnarray}
T^3_{12} &=& -\frac{a\,y''+v\,y'+b\,f\,y'}{b\,y'^2},
\;\;\;\;\;\;\;\;\;\; T^3_{13} = T^3_{14} = T^3_{34}=0,\nonumber \\
T^3_{23} &=&
\frac{2\,y''\,b+a}{b\,y'},\;\;\;\;\;\;\;\;\;\;\;\;\;\;\;\;\;\;\;\;\;\;\;\;\;\;\;
T^3_{24} = \frac{1}{b\,y'},
\end{eqnarray}
and the rest are zero. Since the elements of the group are in the
last row, the coefficients are essential; and thus they could be
absorbed. In the same manner, we can eliminate three other
parameters, so we have,\be
v=-\frac{a\,y''}{y'}-b\,f,\;\;\;\;a=-2\,y''\,b,\;\;\;\;b=-\frac{1}{y'}.\ee
Iterating the procedure, so we have
\begin{eqnarray}
d\omega^1 &=& 0,\;\;\;\;\;\;\;\;\;\;\;\;\;\;d\omega^2 = \omega^2
\wedge\omega^3\nonumber, \\
d\omega^3 &=& \omega^2\wedge\omega^4,\;\;\;\;d\omega^4 =
\sum_{i<j}T^4_{ij}\,\omega^i\wedge\omega^j.
\end{eqnarray}
It is clear that the invariants of this problem are non-zero
coefficients on the fourth line of (34), in other words;
\begin{eqnarray}
I_1 &:=& T^4_{12}=-\frac{1}{y'}\,f_x,\nonumber \\
I_2 &:=&
T^4_{23}=\frac{1}{y'}\,(-3\,f+y'\,f_{y'}+2\,y''\,f_{y''}) \\
I_3 &:=& T^4_{24}=\frac{1}{y'}\,(-3\,y''+y'\,f_{y''}).\nonumber
\end{eqnarray}
\paragraph{Theorem}
{\it A necessary condition that the equations (3) are equivalent
under the time-fixed transformations is that, there exist a
time-fixed transformation $(X,Y)=\Phi(x,y)=(x,\varphi(y))$ such
that $I_1(f)\circ\Phi^{(1)}=I_1(F)$,
$I_2(f)\circ\Phi^{(2)}=I_2(F)$ and $I_3(f)\circ\Phi^{(2)}=I_3(F)$,
where $\Phi^{(i)}:=J^{(i)}\Phi$ is the $i-$jet prolongation of
$\Phi$; in another words
\begin{eqnarray}
\frac{f_x}{y'}\circ\Phi^{(1)}&=&\frac{F_X}{Y'},\nonumber\\
\frac{-3f+y'\,f_{y'}+2y''\,f_{y''}}{y'}\circ\Phi^{(2)}&=&\frac{-3\,F+Y'\,F_{Y'}+2Y''F_{Y''}}{Y'},\\
\frac{-3y'+y'\,f_{y''}}{y'}\circ\Phi^{(2)}&=&\frac{-3Y''+Y'F_{Y''}}{Y'}.
\nonumber\end{eqnarray} }
\section{Sufficient condition}
Achieving the sufficient condition, we use the theory of {\it
$\{e\}$-structures.} Let us ${\cal F}_0:={\cal S}\{I_1,I_2,I_3\}$
be the set of all functions which made by $I_1$ , $I_2$ and
$I_3$. We denote its rank by $k_0$. Since the coframe
$\{\omega^1,\omega^2,\omega^3\}$ is invariant, the derivatives
with respect to them are also invariants. So if $I\in{\cal F}_0$,
we define \be dI = \frac{\partial I}{\partial \omega^1}\,\omega^1
+ \frac{\partial I}{\partial \omega^2}\,\omega^2 + \frac{\partial
I}{\partial \omega^3}\,\omega^3 + \frac{\partial I}{\partial
\omega^4}\,\omega^4\ee then all $\partial I/\partial \omega^i$ are
also invariants. Where,
\begin{eqnarray}
\begin{array}{rclcrcl}
\di \frac{\partial }{\partial \omega^1} &=& \frac{\partial
}{\partial x}, && \di \frac{\partial }{\partial \omega^2} &=& \di
y'\,\frac{\partial }{\partial y}+ y''\,\frac{\partial}{\partial
y'}+f\,\frac{\partial}{\partial y''}, \\[4mm]
\di \frac{\partial}{\partial \omega^3} &=& \di
y'\,\frac{\partial}{\partial y'}+2\,y''\,\frac{\partial}{\partial
y'' }, && \di \frac{\partial}{\partial \omega^4} &=&
y'\,\frac{\partial}{\partial y''} \nonumber
\end{array}
\end{eqnarray}
Now we define,\be {\cal F}_1 = {\cal S}\left\{ I_1,I_2,I_{3};
\frac{\partial I_1}{\partial \omega^1},\frac{\partial
I_1}{\partial \omega^2},\frac{\partial I_1}{\partial
\omega^3},\frac{\partial I_1}{\partial \omega^4} ,
...,\frac{\partial I_3}{\partial \omega^4} \right\},\ee which is
the set of 15 certain functions. Iterating this procedure, we
achieve ${\cal F}_{i}$ and $k_{i}$, for $i=2,3,...$.

By the theory of $\{e\}-$structures, if  $k_i=k_{i+1}$ for some
$i$, then $k_s=k_i$ for all $s\geq i$, moreover $k_i\leq 4$. The
{\it order} of $\{e\}-$structure is the smallest $i$ which
$k_i=k_{i+1}$, and denoted by {\it o}, the value of $k_{i}$ is
also denoted by {\it r} and called the {\it rank} of $\{e\}-$structure.\\
\paragraph{Theorem 4.}
{\it Let ${\cal E}\,:\,y'''=f(x,y,y',y'')$ and $\tilde{\cal
E}\,:\,Y'''=F(X,Y,Y',Y'')$ are two given 3rd order ODEs. Compute
${\cal F}_i$ sets and $k_i$ numbers corresponding to those
equations. The necessary and sufficient condition that these two
equations are equivalent respect to time-fixed transformations is
$(\tilde{x},\tilde{y})=\Phi(x,y)=(x,\varphi(y))$, $\tilde{o}=o$,
$r=\tilde{r}$ and $\tilde{\cal F}_{o+1}={\cal F}_{o+1}\circ
\Phi$.} (See \cite{3}, pp. 271)

\medskip As a result, we have
\paragraph{Conclusion.}
{\it Let $I_1=\alpha$, $I_2=\beta$, and $I_3=\gamma$ are
constants. Then $\gamma=0$, and the ODE is in the form
\begin{eqnarray}
{\cal
E}_{\beta}\;:\;y'''=\frac{3}{2}\,y'\,y''^{2}+y'^3\,h(y)+\frac{
\beta}{2}\,y',
\end{eqnarray}
where $h$ is an arbitrary function of $y$. The ${\cal E}_{\beta}$
is equivalent to ${\cal E}_{\beta'}$ if and only if
$\beta=\beta'$.}

\vspace{2cm}
\noindent Mehdi Nadjafikhah\\
Department of Mathematics,\\
Iran University of Science and Technology,
Narmak-16, Tehran, Iran.\\ E-mail: \verb"m\_nadjafikhah@iust.ac.ir"\\

\noindent Ahmad Reza Forough\\
Department of Mathematics,\\
Iran University of Science and Technology, Narmak-16, Tehran,
Iran.\\ E-mail: \verb"a\_forough@iust.ac.ir"

\begin{thebibliography}{9}
\bibitem{1} {\rm E. Cartan}, {\em Les problemes d'equivalence, oeuvres completes de Elie Cartan}, Vol. III, Center National
de la Recherche Scientifique Paris(1984), pp.1311-1334.
\bibitem{2} {\rm R. Gardner}, {\em The method of Equivalence and it's application Society for Industerial and applied Math.}
, Philadelphia, Pennsylvania, 1989.
\bibitem{4} {\rm M. Nadjafikhah} and {\sc A.R. Forough}, {\em Time fixed geometry of 2nd order
ODEs}, (to appear).
\bibitem{3} {\rm P. Olver}, {\em Equivalence, Invariants and Symmetry}, Cambridge University Press, pp. 252-465.
\end{thebibliography}
\end{document}